\documentclass{amsart}
\usepackage{amsmath}
\usepackage{amssymb}
\usepackage{graphicx}
\usepackage[all,cmtip]{xy}
\input xy
\xyoption{all}
\newtheorem{Lemma}{Lemma}[section]
\newtheorem{Th}[Lemma]{Theorem}

\newenvironment{Proof}{{\sc Proof.}\ }{~\rule{1ex}{1ex}\vspace{0.2truecm}}

\newcommand{\End}{\mbox{\rm End}}

\newcommand{\+}{\oplus}

\begin{document}

    \title[The endomorphism ring of a square-free  injective module]{The endomorphism ring of a square-free \\injective module}

\dedicatory{Dedicated to Professor Alberto Facchini for his 60th birthday}
\author{Mai Hoang Bien}
\address{Mathematisch Instituut, Leiden Universiteit, Niels Bohrweg 1,2333 CA Leiden,The Netherlands.}
\curraddr{Dipartimento di Matematica, Universit\`{a} degli Studi di Padova, Via Trieste 63, 35121 Padova, Italy.}

\email{maihoangbien012@yahoo.com}

      \begin{abstract} Let $M_R$ be an injective right module over a ring $R$. The goal of this paper to prove that the endomorphism ring $S=\End(M_R)$ of $M_R$ is quasi-duo if and only if $M_R$ is square-free. \end{abstract}
    \maketitle
\section{Introduction}\label{1}
In this paper, rings are meant to be associative with identity and modules are understood to be right modules. It is well known that an injective module is indecomposable if and only if its endomorphism ring is a local ring \cite{Mit}. The aim of this paper is to generalize this result from indecomposable modules to square-free ones.

 Recall that an injective $R$-module $M_R$ is said to be {\it square-free} if it contains no direct summand isomorphic to $X\+X$ for some direct summand $0\ne X$ of $M_R$. A ring $S$ is called {\it left quasi-duo} (resp. {\it right quasi-duo}) if every maximal left ideal (resp. maximal right ideal) of $S$ is a two-sided ideal.  If $S$ is both left and right quasi-duo then $S$ is said to be {\it quasi-duo}. The class of quasi-duo rings is not ``small". For example, commutative rings and local rings are quasi-duo. Factors of a quasi-duo ring are also quasi-duo. In \cite{Yu}, the author proved that $S$ is quasi-duo if and only if $S/J(S)$ is quasi-duo. Here, $J(S)$ is the Jacobson radical of $S$. The direct product of a family of rings is quasi-duo if and only if each ring is quasi-duo. The matrices ring $M_n(D)$ of degree $n$ over a division ring $D$ is quasi-duo if and only if $n=1$. Hence, a direct product of finitely many division rings is quasi-duo. Therefore, rings of finite type are quasi-duo \cite{Fac3}. In particular, the endomorphism rings of univerial modules, couniform projective modules, cyclically presented modules over a local ring, and the kernels of morphisms between indecomposable injective modules  are quasi-duo (see \cite{Fac3}). Examples of right quasi-duo rings which are not left quasi-duo are unknown \cite[Question 7.7]{Lam2}. Of course, if there exists such ring $R$ then the opositive ring $R^{op}$ of $R$ is left quasi-duo and not right quasi-duo.
 
 Let $M_R$ be an injective module over $R$. In this paper, we will show that $M_R$ is square-free if the endomorphism ring $S: =\End(M_R)$ of $M_R$ is either left quasi-duo or right quasi-duo (Theorem~\ref{2.1}). Conversely, if $M_R$ is square-free then $S$ is quasi-duo. In fact, we prove that if $M_R$ is square-free then every one-sided ideal of $S$ containing $J(S)$ is two-sided (Theorem~\ref{2.6} and~\ref{2.8}).  

The symbols and notations we use in this paper are familiar to module and ring theory. $\subseteq _e $ will denote an essential submodule, and $E(N)$ will denote an injective envelope of a module $N$. For any direct summand $N$ of a $R$-module $M_R$, $\iota_N$ and $\pi_N$, respectively, are the embedding of $N$ into $M_R$ and the projection of $M_R$ onto $N$. 

\section{Results}\label{2}
\begin{Th}\label{2.1} Let $M_R$ be an injective $R$-module. If the endomorphism ring $S=\End(M_R)$ of $M_R$ is either left quasi-duo or right quasi-duo then $M_R$ is square-free.
\end{Th}

\begin{Proof} Assume that $S$ is left quasi-duo and there exist direct summands $A,B,C$ of $M_R$ such that $A\+B\+C=M_R$ and $A$ is isomorphic to $B$. Call $I$ a maximal left ideal of $S$ containing $s=\iota_{B\+C}\pi_{B\+C}$. Then $I$ is a two-sided ideal of $S$ by hypothesis. Let $\alpha$ be an isomorphism from $A$ to $B$. Define $f:=\iota_B\alpha\pi_A$ and $g:=\iota_A\alpha^{-1}\pi_B+\iota_C\pi_C$. Then $s f+g s\in I$ and 

$$\ker (s f+gs) =\{\, a+b+c\in A\+B\+C\mid (s f+gs)(a+b+c)=0 \,\}$$ $$=\{\, a+b+c\in A\+B\+C \mid \alpha (a)+\alpha^{-1}(b)+c=0 \,\}=0.$$ Hence, $s f+gs$ is injective. Consider the diagram 

$$\xymatrixcolsep{5pc}
\xymatrix{
0 \ar[r] &{A\+B\+C=M_R}\ar[r]^{sf+gs}  \ar[d]^{1} &M_R \ar@{.>}[ld]^h     \\
&M_R}
$$Since $M_R$ is injective, there exists morphism $h$ of $S$ such that $h(s f+gs)=1$, which imlies $I=S$. Contradiction. Therefore, if $S$ is left quasi-duo then $M_R$ is square-free.

Similarly to the right side. 
 \end{Proof}

\begin{Lemma}\label{2.2} Let $M_R$ be a  square-free injective $R$-module. If $N_1, N_2$ are two isomorphic direct summands of $M_R$ then $N_1$ and $N_2$ are two injective envelopes of $N_1\cap N_2$.
\end{Lemma}
\begin{Proof} Let $A_1, A_2$ be injective envelopes of $N_1\cap N_2$ respectively in $N_1$, $N_2$ and $B_1,B_2$ be respectively direct summands of $N_1,N_2$ such that $$N_1=A_1\+B_1, N_2=A_2\+B_2.$$ Since $N_1\cong N_2$ and $B_1\cap A_2=0$,   $B_1$ is isomorphic to a direct summand $C$ of $N_2$ and since $B_1\cap N_2=0$, $B_1+N_2=B_1\+N_2$ is injective, which implies that $M_R$ contains $B_1\+C$, with $B_1\cong C$, as a direct summand. Therefore $B_1=0$. Similarly, $B_2=0$. Thus, $N_1, N_2$ are two injective envelopes of $N_1\cap N_2$.\end{Proof}

\begin{Lemma}\label{2.3} Let $M_R$ be a square-free injective $R$-module with $M_R=M_1\+M_2$. If $N$ is a direct summand of $M_R$ then $N$ is an injective envelope of $(N\cap M_1)\+(N\cap M_2)$.
\end{Lemma}
\begin{Proof} Let $A$ be an injective envelope of $N\cap M_1$ in $N$ and $B$ be a direct summand of $N$ such that $N=A\+B$. Consider $\pi:=\pi_{M_2}$, the projection of $M_1\+M_2$ onto $ M_2$, and the morphism $\pi_{\left.   \right| B}:B\to M_2$ restricted on $B$ of $\pi$. Because $B\cap M_1=0$, ${\left. \pi  \right|_B}$ is injective. Hence, $B$ is isomorphic to a direct summand $C$ of $M_2$. By Lemma~\ref{2.2}, $B$ is an injective envelope of $B\cap C$. Therefore, $B$ is an injective envelope of $N\cap M_2=B\cap M_2\supseteq B\cap C$. This implies $N$ is an injective envelope of $(N\cap M_1)\+(N\cap M_2)$. 
\end{Proof}

\begin{Lemma}\label{2.4} Let $M_R$ be an injective $R$-module, $S=\End(M_R)$ be the endomorphism ring of $M_R$ and $J(S)$ be the Jacobson radical of $S$. For any element $f\in S$, there exist $e_1,e_2, g_1,g_2, h_1,h_2\in S$ and $ i_1, i_2, j_1, j_2\in J(S)$ such that $e_1,e_2$ are idempotents and $$e_1=fg_1+i_1, f=e_1h_1+j_1,$$ $$e_2=g_2f+i_2, f=h_2e_2+j_2.$$
\end{Lemma}
\begin{Proof} It is well-known that $S/J(S)$ is a Von Neumann regular ring and any idempotent of $S/J(S)$ can be lifted to an idempotent of $S$ (see \cite[Theorem 13.1]{Lam}). Hence, for any $f\in S$, there exists idempotent $e_1\in S$ such that $$(fS+J(S))/J(S)=(e_1S+J(S))/J(S).$$ Therefore, $e_1=fg_1+i_1, f=e_1h_1+j_1$ for some $g_1,h_1\in S$ and $i_1,j_1\in J(S)$. Similarly to $e_2=g_2f+i_2, f=h_2e_2+j_2.$   
\end{Proof}

\begin{Lemma}\label{2.5} Let $M_R$ be an injective $R$-module, $S=\End(M_R)$ be the endomorphism ring of $M_R$ and $e$ be an idempotent of $S$. For any element $f\in S$, $f\in eS$ if and only if $f(M_R)\subseteq e(M_R)$. 
\end{Lemma}
\begin{Proof} Assume that $f=eh$ for some $h\in S$. Then $f(M_R)=eh(M_R)\subseteq e(M_R)$. Conversely, assume that $f(M_R)\subseteq e(M_R)$. Then for any $x\in M_R$, $f(x)=e(y)$ for some $y\in M_R$. Hence, $$f(x)=e(y)=e(e(y))=e(f(x)).$$ Therefore, $f=ef\in eS$.
\end{Proof}

\begin{Th}\label{2.6} Let $M_R$ be an injective $R$-module and $S=\End(M_R)$ be the endomorphism ring of $M_R$. If $M_R$ is square-free then every right ideal of $S$ containing $J(S)$ of $S$ is a two-sided ideal. In particular, if $M_R$ is square-free then $S$ is right quasi-duo. 
\end{Th}

\begin{Proof} Let $I$ be a right ideal of $S$ containing $J(S)$. Let $f\in I$ and $\phi\in S$. We must show that $\phi f\in I$. By Lemma~\ref{2.4}, there exist $e,g,h\in S$ and $i,j\in J(S)$ such that $e=fg+i$ is an idempotent and $f=eh+j$. We have $e$ belongs to $I$ and $\phi f=\phi e h+\phi j$, so that it suffices to show that $\phi e\in I$. Indeed, let $N'$ be a direct summand of $M_R$ such that $M_R=N\+N'$ with $N=e(M_R)$. If set $M_1$ is an injective envelope of $\ker \phi$, then there exists a direct summand $M_2$ of $M_R$ such that $M_R=M_1\+M_2$.  Let $N_1, N_2$ be respectively injective envelopes of $N\cap M_1, N\cap M_2$ in $N$. By Lemma~\ref{2.3}, we may assume that $$M_R=N\+N'=N_1\+N_2\+N'.$$  Consider the morphism $\phi _{\left. \right|{N_2}}:N_2\to N_1\+N_2\+N'$ restricted on $N_2$ of $\phi$. It is easy to check that $\phi _{\left. \right|{N_2}}$ is injective.  Hence, $N_2\cong \phi _{\left. \right|{N_2}}(N_2)$ and by Lemma~\ref{2.2}, $A=N_2\cap \phi _{\left. \right|{N_2}}(N_2)$ is an essential submodule of $N_2$. Set $B=\phi _{\left. \right|{A}}^{-1}(A)$ with $\left. \phi\right|_{A}:A\to N_2$ is the morphism restricted on $A$ of $\phi$. Then $B$ is also an essential submodule of $N_2$. Write $\psi: N_2\to N_2$ for a morphism extending of $\phi _{\left. \right|{A}}$ and let $\psi'=\iota_{N_2}\psi\pi_{N_2}$. One has $\psi'(M_R)\subseteq N_2\subseteq N=e(M_R)$, which implies from Lemma~\ref{2.5} that $\psi'\in eS\subseteq I$. Moreover, since for any $a+b+c\in (\ker \phi\cap N_1)\+B\+N'$, $(\phi e-\psi')(a+b+c)=0$ and $(\ker \phi\cap N_1)\+B\+N'$ is essential in $M_R$, $\phi e-\psi'$ belongs to $J(S)$. This shows that $\phi e\in I$.
\end{Proof}

\begin{Lemma}\label{2.7} Let $M_R$ be an injective $R$-module, $S=\End(M_R)$ be the endomorphism ring of $M_R$. Then for any  two elements  $f,g$ of $S$,  $g=hf+j$ for some $h\in S, j\in J(R)$ if and only if $\ker f\cap A \subseteq \ker g$ for some essential submodule $A$ of $M_R$. 
\end{Lemma}
\begin{Proof}Assume that $g=hf+j$ for some $h\in S$ and $j\in J(S)$. Then $\ker g\supseteq \ker (hf)\cap \ker j \supseteq \ker f\cap A$ with $A=\ker(j)\subseteq_e M_R$.

Conversely, assume that $\ker f\cap A\subseteq \ker g$ for a submodule $A\subseteq _eM_R$. Let $M_1$ be an injective envelope of $\ker g$ in $M_R$, and $N_1$ be an injective envelope of $\ker f \cap A$ in $M_1$. Assume that $M_R=N_1\+N_2\+M_2$ for some direct summands $N_2, M_2$ respectively of $M_1$ and $M_R$. Consider the diagram

$$\xymatrixcolsep{5pc}
\xymatrix{
0 \ar[r] &{N_2\+M_2}\ar[r]^{f_{\left.  \right|{N_2 \oplus M_2}}}  \ar[d]^{g_{\left.  \right|{N_2 \oplus M_2}}} &M_R \ar@{.>}[ld]^h     \\
&M_R}
$$ Here, ${f_{\left.  \right|{N_2 \oplus M_2}}}, {g_{\left.  \right|{N_2 \oplus M_2}}}$ are the morphisms respectively restricted on $N_2\+M_2$  of $f,g$. Since $M_R$ is injective, there exists $h:M_R\to M_R$ such that $${g_{\left.  \right|{N_2 \oplus M_2}}}=h.{f_{\left.  \right|{N_2 \oplus M_2}}}.$$ Because $(g-hf)(a+b+c)=0$ for any  $a+b+c\in (\ker f\cap A)\+N_2\+M_2$ which is essential in $M_R$, $j=g-hf\in J(S)$.
\end{Proof}

\begin{Th}\label{2.8} Let $M_R$ be an injective $R$-module and $S=\End(M_R)$ be the endomorphism ring of $M_R$. If $M_R$ is square-free then every left ideal of $S$ containing the Jacobson radical $J(S)$ of $S$ is a two-sided ideal. In particular, if $M_R$ is square-free then $S$ is left quasi-duo. 
\end{Th}
\begin{Proof}Let $I$ be a left ideal of $S$ containing $J(S)$. Let $f,\phi$ be elements of $ S$ with $f\in I$. We must show that $f\phi \in I$. By Lemma~\ref{2.4}, there exist $e,g,h\in S$ and $i,j\in J(S)$ such that $e=gf+i$ is an idempotent and $f=he+j$. Hence, $e\in I$ and $f\phi =he \phi+j\phi$. It suffices to show that $e\phi \in I$. Indeed, set $N_1=\ker e$ and let $N_2$ be a direct summand of $M_R$ such that $M_R=N_1\+N_2$. Consider $\psi:=\pi_{N_2}\phi\iota_{N_1}:N_1\to N_2$. If set $N_1'$ is an injective envelope of $A:=\ker \phi$ in $N_1$ then exists a direct summand $N_1''$ of $M_R$ such that $N_1=N_1'\+N_1''$. Since $\ker \psi_{\left.  \right|{{N_1''}}}=N_1''\cap A=0$, $\psi_{\left.  \right|{{N_1''}}}$ is injective. This implies $N_1''$ is isomorphic to a direct summand $C$ of $N_2$. Hence, $M_R$ contains a direct summand isomorphic to $C\+C$. By hypothesis, $C=0$. In other words, $A$ is essential in $N_1$. Now, one has that $\phi$ can be written as the matrix $\phi=\left( {\begin{array}{*{20}{c}}
{\pi_{N_1}\phi\iota_{N_1}}&{\pi_{N_1}\phi\iota_{N_2}}\\
{\pi_{N_2}\phi\iota_1}&{\pi_{N_2}\phi\iota_{N_2}}
\end{array}} \right)$. Then $\ker (e\phi)=\phi^{-1}(\ker e)=\phi^{-1}(N_1)\supseteq \ker(\pi_{N_2}\phi\iota_1)\+\ker(\pi_{N_2}\phi\iota_{N_2})\supseteq A=(A\+N_2)\cap N_1=(A\+N_2)\cap \ker e$. Since $A\+N_2$ is essential in $M_R$, by Lemma~\ref{2.7}, there exist $h\in S, j\in J(S)$ such that $e\phi=he+j$. Thus $e\phi\in I$.\end{Proof}

\end{document}